\documentclass[12pt]{article}
\begin{document}
\newcommand{\qed}{\ \
\mbox{\rule{8pt}{8pt}}\vspace{0.3cm}\newline}
\newcommand{\ia}{{\bf I}_{A_i}}
\newcommand{\ba}{\widetilde{{\bf KC}}_i}
\newcommand{\bba}{{\bf KC}_i}
\newcommand{\ra}{\longrightarrow}
\newcommand{\pe}{{\cal P}}
\newcommand{\der}{{\cal DP}_d}
\newcommand{\adp}{{\cal DP}_d^{af}}
\newcommand{\ot}{\otimes}
\newcommand{\rec}{\raisebox{-1ex}{\ $\stackrel{\textstyle{\stackrel{\textstyle{\longleftarrow}}{\longrightarrow}}}{\longleftarrow}$\ }}
\title{Affine strict polynomial functors and formality}
\author{Marcin Cha\l upnik
\thanks{The author was supported by the  grant (NCN) 2011/01/B/ST1/06184.}\\
\normalsize{Institute of Mathematics, University of Warsaw,}\\
\normalsize{ul.~Banacha 2, 02--097 Warsaw, Poland;}\\
\normalsize{e--mail: {\tt mchal@mimuw.edu.pl}}}
\date{\mbox{}}
\newtheorem{prop}{Proposition}[section]
\newtheorem{cor}[prop]{Corollary}
\newtheorem{theo}[prop]{Theorem}
\newtheorem{lem}[prop]{Lemma}
\newtheorem{defi}[prop]{Definition}
\newtheorem{defipro}[prop]{Definition/Proposition}
\newtheorem{fact}[prop]{Fact}
\newtheorem{exa}[prop]{Example}
\newcommand{\ka}{{\mbox {\bf k}}}
\newcommand{\ca}{\mbox{$\cal{A}$}}
\maketitle
\begin{abstract}
We introduce the notion of affine strict polynomial functor. We show how this concept helps to understand
homological behavior of the operation of Frobenius twist
in the category of strict polynomial functors over a field
of positive characteristic. We also point out for an analogy
between our category and the category of representations  of the group
of algebraic loops on $GL_n$.
\\\mbox{}\vspace{0.1cm}
\end{abstract}
\section{Introduction}
In the present paper we study homological algebra in the category
${\cal P}_d$ of strict polynomial functors of degree $d$ over a field of characteristic $p>0$.
We introduce a new type of strict polynomial functors we call the affine strict polynomial functors. This was motivated by observation that the Ext--groups  of the form
$\mbox{Ext}^*_{{\cal P}_{pd}}(F^{(1)},G)$ where $F^{(1)}$ is the Frobenius twist of $F\in{\cal P}_d$ are equipped with certain extra structure.
It has become more and more apparent that for needs of the program started in [C1], [C3] aiming at computing the Ext groups between
strict polynomial functors important in representation theory,
understanding this extra structure is necessary. This structure comes, roughly speaking, from the fact that the Frobenius twist
functor $I^{(1)}\in{\cal P}_p$ has non--trivial endo--Ext--groups and these Ext--groups act (nonlinearly) on $F^{(1)}:=F\circ I^{(1)}$. However, for technical reasons, it is more convenient to look at the operation adjoint to the twisting. Namely, precomposing with $I^{(1)}$ is an exact operation, hence it gives rise to a functor between (bounded) derived categories:
\[{\bf C}:{\cal DP}_d\ra{\cal DP}_{pd}.\]
It was shown in [C4] that this functor has a right adjoint
\[{\bf K}^r:{\cal DP}_{pd}\ra{\cal DP}_{d}\]
called there the derived right Kan extension. Our idea is to
factorize ${\bf K}^r$ through certain richer triangulated category
${\cal DP}_d^{af}$ which should be thought of as  the derived category of the category ${\cal P}_d^{af}$ of affine strict polynomial functors of degree $d$. Now since
\[\mbox{Ext}^*_{{\cal P}_{pd}}(F^{(1)},G)\simeq
\mbox{HExt}^*_{{\cal P}_{d}}(F,{\bf K}^r(G))\]
and ${\bf K}^r(G)$ comes from this richer category we see where extra
structure comes from. In fact this enriched category ${\cal P}_d^{af}$ has quite transparent interpretation: it is  the category
of strict polynomial functors but from (certain subcategory of) the category of (graded free
 finitely generated) $A$--modules where $A:=\mbox{Ext}^*_{{\cal P}_p}
 (I^{(1)},I^{(1)})$. However, since our construction is technically a bit involved, let me present here some informal review of it.  We first recall from [C4] definition
 of ${\bf K}^r$. Namely we have
 \[{\bf K}^r(F)=\mbox{RHom}_{{\cal P}_{pd}}(\Gamma^d(I^{(1)}\ot I^*),F)\]
 where  $\Gamma^d(I^{(1)}\ot I^*)$ is certain functor in two variables. Thus
 by taking RHom with respect to the first variable we are left with a complex of functors in the second variable. Now we can see some extra structure on $H^*({\bf K}^r(F))$: the Ext--groups
 $\mbox{Ext}^*_{{\cal P}_{pd}}(\Gamma^d(I^{(1)}\ot I^*),\Gamma^d(I^{(1)}\ot I^*))$ act on it.
 These Ext--groups will turn out to be Hom--spaces in the source of our functor category
 ${\cal P}_d^{af}$.
 However, the situation is still not satisfactory since we want these Ext--groups acting on ${\bf K}^r(F)$, not
 merely on its cohomology.
 To resolve this problem we need some sort of formality.
 Let us make this more precise. We fix a projective resolution $X$ of $\Gamma^d(I^{(1)}\ot I^*)$ and introduce an auxiliary
 differential graded (=DG) category (denoted by $\Gamma^d {\cal V}_X$)  whose Hom spaces are complexes
 $\mbox{Hom}_{{\cal P}_{pd}}(X(-,V'),X(-,V)$. This category genuinely  acts on ${\bf K}^r(F)$ and our Ext--groups form the cohomology category of it. Now the main point is Theorem 4.2 which says that the DG category
 $\Gamma^d {\cal V}_X$ is formal. This allows us to endow ${\bf K}^r(F)$ with
 a  structure which a priori existed only on cohomology. \newline
 The article is organized as follows. In Section 2 we introduce the category ${\cal P}_d^{af}$  and establish its basic properties. Section 3
 recalls generalities on derived categories of DG  categories and discusses certain finiteness assumption which is useful in our situation.
 In Section 4 we introduce the category $\Gamma^d {\cal V}_X$ and prove its formality. Finally, in Section 5 we construct
 the affine derived right Kan extension ${\bf K}^{af}$ and relate it to ${\bf K}^r$. Theorem 5.1 which establishes the fundamental properties
  of ${\bf K}^{af}$ is the main result of the paper. In the last section we
 put our work into a wider context and discuss some possible further developments. In particular we observe a formal analogy between ${\cal P}_d^{af}$ and the category of representations of the groups of algebraic loops on $GL_n$.\newline {\em Acknowledgements} I am grateful to Julian K\"ulshammer for turning
 my attention to [Kle].
\section{The category of affine strict polynomial functors}
We fix a field {\bf k} of characteristic $p>0$. Let ${\cal V}$ stands for the category of finite dimensional \ka--spaces and
${\cal V}^{f+}$ be the category
of {\bf Z}--graded bounded below {\bf k}--spaces, finite dimensional in each degree.
 Let $A$ denote the graded \ka-algebra $A:={\bf k}[x]/x^p$ for $x$ of degree $2$.
 By the classical computation [FS, Th.~4.5]: $A\simeq\mbox{Ext}^*_{{\cal P}_p}(I^{(1)},I^{(1)})$.
  We consider the category ${\cal V}_A$
 which is the full subcategory of the category of graded $A$--modules consisting of objects of the form $V\ot A$ for $V\in {\cal V}$.
 Hence
 \[\mbox{Hom}_{{\cal V}_A}(V\ot A,W\ot A)=\mbox{Hom}_A(V\ot A,W\ot A)\simeq\mbox{Hom}(V,W)\ot A\]
(unless otherwise stated all (graded) linear and tensor operation are taken over \ka).\newline
We shall introduce certain version of the Friedlander--Suslin strict polynomial functors which, roughly speaking, correspond to the functors from
$A$--modules to \ka--modules. Technically, it will be convenient to adopt the  approach to the strict polynomial functors
due to Pirashvili  which allows to interpret strict polynomial functors as genuine functors.
We recall (see e.g. [FP, Sect.~3]) that one considers the category $\Gamma^d {\cal V}$ whose objects are  finite dimensional \ka--spaces but
\[\mbox{Hom}_{\Gamma^d {\cal V}_{\bf k}}(V,W):=\Gamma^d(\mbox{Hom}(V,W))\]
where $\Gamma^d(X):=(X^{\ot d})^{\Sigma_d}$ is the space of symmetric d--tensors on a \ka--space $X$. Then it is easy to see that a strict polynomial
functor homogeneous of degree $d$ is nothing but a genuine \ka--linear functor from $\Gamma^d {\cal V}$ to ${\cal V}$.\newline
Now we introduce the category $\Gamma^d {\cal V}_{A}$ with the objects the same as in ${\cal V}_A$ and
 \[\mbox{Hom}_{\Gamma^d {\cal V}_{A}}(V\ot A,W\ot A):=\Gamma^d(\mbox{Hom}(V,W)\ot A).\]
  The categories $\Gamma^d {\cal V}_{A}, {\cal V}^{f+}$ are graded \ka--linear categories (by this we mean that the Hom sets are graded \ka--linear
 spaces and composition preserves this structure). We call functor between graded categories a graded functor if its action on the Hom spaces preserves grading.
 Now we are ready for defining our functor category.
 \begin{defipro}
An affine strict polynomial functor $F$ homogeneous of degree $d$ is a graded \ka--linear functor
\[F:\Gamma^d {\cal V}_{A}\ra {\cal V}^{f+}.\]
The affine strict polynomial functors homogeneous of degree $d$ with morphisms being natural transformations form a graded abelian category ${\cal P}_d^{af}$.
\end{defipro}
Let me at this point comment on our notational conventions. In fact the objects of ${\cal V}$ are indexed  just by vector spaces. Nevertheless, we prefer
to use label  $V\ot A$ instead of $V$. The reason is that when we construct functors on  ${\cal V}_A$ we usually somehow use $A$--structure on Hom--spaces. Even if not, like in the forgetting functor $I\in{\cal P}_1^{af}$ which sends $V\ot A$ to itself, the apparently simpler notation would
produce quite strange formula: $I(V):=V\ot A$. Hence forgetting would rather look like inducing which would be quite confusing.\newline
Like in any functor category, for any $U\in{\cal V}$ we have the representable functor $h^{U\ot A}\in {\cal P}_d^{af}$ given by the formula
\[V\ot A\mapsto \mbox{Hom}_{\Gamma^d {\cal V}_{A}}(U\ot A,V\ot A)\]
and the co--representable functor $c_{U\ot A}^*\in {\cal P}_d^{af}$ given by the formula
\[V\ot A\mapsto \mbox{Hom}_{\Gamma^d {\cal V}_{A}}(V\ot A,U\ot A)^*\]
where $(-)^*$ stands for the graded \ka--linear dual.
Now we have
\begin{prop}
\begin{enumerate}
\item
There are  natural in $U\ot A$ isomorphisms
\[\mbox{Hom}_{{\cal P}_d^{af}}(h^{U\ot A},F)\simeq F(U\ot A)\]
\[\mbox{Hom}_{{\cal P}_d^{af}}(F,c_{U\ot A}^*)\simeq (F(U\ot A))^*\]
for any $F\in{\cal P}_d^{af}$.\newline
\item
Moreover, the map $\Theta: h^{U\ot A}\otimes F(U\ot A)\ra F$ adjoint to the map $F_{U\ot A,V\ot A}$ giving the action of $F$ on morphisms is surjective, provided that $\dim(U)\geq d$.
\item
If $\dim(U)\geq d$ then  $h^{U\ot A}$ is a projective generator of ${\cal P}_d^{af}$, $c_{U\ot A}^*$ is an injective generator of ${\cal P}_d^{af}$
\end{enumerate}
\end{prop}
{\bf Proof: }The first part is just the Yoneda lemma.\newline
The surjectivity of $\Theta$ is proved analogously to [FS, Th. 2.10].\newline
The projectivity of $h^{U}$ follows from part 1, the fact that it is a generator follows from part 2. The statements about $c^*_{U\ot A}$ are proved analogously.\qed
Since ${\cal V}_A$ is sort of scalar extension of ${\cal V}$, one can  expect some inducing/forgetting adjunction between our functor categories. Indeed, we have the
functors
$z:{\cal V}_A\ra {\cal V}$, $t:{\cal V}\ra {\cal V}_A$ given on objects by formulae: $z(V\ot A):=V\ot A$, $t(V):=V\ot A$. The action on Hom--spaces is the following:
\[z_{V,W}: \mbox{Hom}(V,W)\ra \mbox{Hom}(V,W)\ot A\]
is the natural embedding;
\[t_{V,W}: \mbox{Hom}(V,W)\ot A\simeq \mbox{Hom}_A(V\ot A,W\ot A)\ra\mbox{Hom}(V\ot A,W\ot A)\]
is the embedding again. These functors are nothing but forgetting and induction functors.  These functors extend naturally to
$z:\Gamma^d{\cal V}_A\ra \Gamma^d{\cal V}$, $t:\Gamma^d{\cal V}\ra \Gamma^d{\cal V}_A$ which will be denoted by the same letters $z,t$. Now we shall
show that that the assigning $F\mapsto F\circ z$ gives the functor $z^*:{\cal P}_d^{f+}\ra {\cal P}_d^{af}$ where ${\cal P}_d^{f+}$
is the category of graded bounded below strict polynomial functors of degree $d$. This is not entirely obvious as can  already be
seen for $F$ concentrated in degree $0$. Then taking $z^*(F)(V):=F(V\ot A)$ as concentrated in degree $0$ does not produce a graded functor. The correct
approach relies on the fact (specific to strict polynomial functors), that $F$ can be naturally extended to the functor
$F^{gr}: \Gamma^d{\cal V}^{f+}\ra {\cal V}^{f+}$ (see e.g. [T2, Sect.~2.5]). We recall that for the graded space $V=\bigoplus V^j$ we have decomposition
$F(V)=\bigoplus F_{\gamma}$ where for $\gamma=(\gamma_1,\ldots,\gamma_s)$ with $\sum \gamma_j=d$ $F_{\gamma}(\bigoplus V^j)$ is the sub--s--functor of
$F(\bigoplus V^j)$ of degree $\gamma_j$ in $V^j$. Then we put to $F_{\gamma}$ degree $\sum j\gamma_j$. In general, for $F=\bigoplus F^s\in {\cal P}_d^{af}$
we assign to $F^s(V\ot A)$ degree $s+\sum_j j\gamma_j$. Now it is easy to see that $z^*$ is well defined. We analogously define $t^*$ as preocmposition with $t$.
Now, as one can expect, we have
\begin{prop}
\begin{enumerate}
\item
$z^*$ preserves representable objects i.e. \[z^*(\Gamma^{d,U})=h^{U\ot A}.\]
where $\Gamma^{d,U}\in{\cal P}_d$ is defined as $V\mapsto\mbox{Hom}_{\Gamma^d{\cal V}}(U,V)=\Gamma^d(U^*\ot V)$.
\item
The functor $t^*$ is right adjoint to $z^*$.
\end{enumerate}
\end{prop}
{\bf Proof: } We have
\[z^*(\Gamma^{d,U})(V\ot A)=\Gamma^{d,U}(V\ot A)=\Gamma^d(\mbox{Hom}(U,V\ot A))\simeq\Gamma^d(\mbox{Hom}(U,V)\ot A)\simeq h^{U\ot A}(V\ot A),\]
thus getting the first part. Since $\Gamma^{d,U}$ generate ${\cal P}_d$, in order to get the second part, it suffices to obtain a
natural isomorphism
\[\mbox{Hom}_{{\cal P}_d}(\Gamma^{d,U},t^*(F))\simeq \mbox{Hom}_{{\cal P}_d^{af}}(z^*(\Gamma^{d,U}),F).\]
for any $F\in{\cal P}_d^{af}$.
Now by the Yoneda lemma we have
\[\mbox{Hom}_{{\cal P}_d}(\Gamma^{d,U},t^*(F))\simeq t^*(F)(U)=F(U\ot A).\]
By using the first part we get
\[\mbox{Hom}_{{\cal P}_d^{af}}(z^*(\Gamma^{d,U}),F)\simeq \mbox{Hom}_{{\cal P}_d^{af}}(h^{U\ot A},F)\simeq F(U\ot A).\]
\qed
The functor $z^*$ provides  a lot of examples of affine strict polynomial functors. We shall occasionally denote
$z^*(F)$ as $F^{af}$.
E.g. we have affine versions of tensor functors: $(I^d)^{af}(V\ot A):=(V\ot A)^{\ot d}$, $(S^d)^{af}(V\ot A):=S^d(V\ot A)$ etc. Perhaps more interesting
are objects in ${\cal P}_d^{af}$ which do not come from ${\cal P}_d$. The most fundamental examples are
$\chi_j\in {\cal P}_1^{af}$ for $0\leq j\leq p-1$ given by \[\chi_j(V\ot A):=(x^{j}.\ (V\ot A))/(x^{j+1}.\ (V\ot A))\simeq V[-2j].\]
We have also the ``affine Kuhn duality'' in ${\cal P}_d^{af}$, although  its construction requires some explanation. The reason is
that the formula $F^{\#}(V\ot A):=(F((V\ot A)^*))^*$ technically does not make sense since $(V\ot A)^*$ is not an object of ${\cal V}_A$. However,
since $A^*\simeq A[2(p-1)]$ as graded $A$--modules, we just have $(V\ot A)^*\simeq V^*\ot A^*\simeq V^*\ot A[2(p-1)]$. Thus we formally define
\[F^{\#}(V\ot A):=F^{gr}(V^*\ot A[2(p-1)])^*\]
where $F^{gr}$ is the extension of $F$ to the graded spaces we have discussed earlier. Also later we will sometimes apply our functors to graded
spaces tacitly assuming using their graded extensions.
 As one can expect, the affine Kuhn duality takes representable
functors to co--representable ones. However, also here, some shifting phenomena emerge. The best way to capture them is to allow representable and
co--representable functors to be labeled by graded spaces, which again is justified by using graded extensions of functors. Now we gather the basic properties
of the Kuhn duality
\begin{prop}\mbox{}
\begin{enumerate}
\item
$(h^{U\ot A})^{\#}\simeq c_{U^*\ot A^*}^*\simeq c_{U^*\ot A}^*[-2(p-1)d]$
\item
$(\chi_j)^{\#}=\chi_{p-1-j}$
\item
For any $F\in{\cal P}_d$, $z^*(F^{\#})\simeq (z^*(F))^{\#}$.
\item
The functor $h^*:=(-)^{\#}\circ t^*\circ (-)^{\#}$ is left adjoint to $z^*$.
\item
Explicitly: $h^*(F)(V)\simeq F(V\ot A^*)$, hence $h^*\simeq t^*[2(p-1)d]$.
\end{enumerate}
\end{prop}
{\bf Proof: } For the first part we compute
\[(h^{U\ot A})^{\#}(V\ot A)=(\Gamma^d(\mbox{Hom}_A(U\ot A,(V\ot A)^*)))^*\simeq (\Gamma^d(\mbox{Hom}_A(V\ot A,(U\ot A)^*)))=\]
\[c_{U^*\ot A^*}(V\ot A).\]
The isomorphism $c_{U^*\ot A^*}^*\simeq c_{U^*\ot A}^*[-2(p-1)d]$ follows from properties of graded extension.\newline
The second part follows from the fact that $\chi_j(V\ot A)\simeq V[-2j]$.
The third part is obvious, the fourth part follows formally from the third part part and $\{z^*,t^*\}$ adjunction. The fifth part is obvious.\qed
We finish reviewing basic properties of the category ${\cal P}_d^{af}$ by investigating properties of the functor $ev_n$ sending $F$ to $F(A^n)\simeq F({\bf k}^n\ot A)$. Let $S_{d,n}^{af}$ denote the graded \ka--algebra
$\Gamma^d(\mbox{End}_A(A^n))\simeq\Gamma^d(\mbox{End}({\bf k}^n)\ot A)$.
Then by Prop.~2.2.1 $\mbox{Hom}_{{\cal P}_d^{af}}(h^{A^n},h^{A^n})\simeq S_{d,n}^{af}$. Therefore, since by Prop.~2.2.1 again $F(A^n)\simeq \mbox{Hom}_{{\cal P}_d^{af}}(h^{A^n},F)$,
$F(A^n)$ is endowed with a natural structure of graded $S_{d,n}^{af}$--module.
\begin{prop}
If $n\geq d$ then
\[ev_n: {\cal P}_d^{af}\ra S_{d,n}^{af}\mbox{-mod}^{f+},\]
where $S_{d,n}^{af}\mbox{-mod}^{f+}$ is the category of bounded below finite dimensional in each degree graded $S_{d,n}^{af}$--modules,
is an equivalence of graded abelian categories.
\end{prop}
{\bf Proof: } is analogous to that of [FS, Th.~3.2]. We assign to $M\in S_{d,n}^{af}\mbox{-mod}^{f+}$ the affine functor $\Phi(M)$ given by the formula
\[V\ot A\mapsto \Gamma^d(Hom_{A}(A^n,V\ot A)\otimes_{S_{d,n}^{af}} M.\]
It is easy to see that $ev_n\circ \Phi\simeq Id_{S_{d,n}^{af}\mbox{-mod}_{gr}^{f+}}$ and that $\Phi\circ ev_n(h^{A^n})\simeq h^{A^n}$. This, since
$h^{A^n}$ generates ${\cal P}_d^{af}$, shows that $\Phi$ is quasi--inverse of $ev_n$.\qed
\section{Deriving ${\cal P}_d^{af}$}
In this section we introduce the derived category of ${\cal P}_d^{af}$. For ${\cal P}_d^{af}$ is graded abelian category, it is natural to look at
it as a DG category (with trivial differentials) and use the formalism of derived categories of DG categories. Since DG homological algebra is not
as well known as its abelian counterpart we start with recalling
some standard constructions concerning DG functor categories and their derived categories.  Our main reference in this section is a classical paper [K1], in particular we borrow notation from there.\newline
By DG category we mean  a \ka--linear category whose Hom--sets are naturally (possibly unbounded) cohomological complexes. Thus a DG category is in particular a graded category. The simplest example of DG category is the category $\mbox{Dif}(\ka)$ of complexes of \ka--vector spaces with internal Hom complexes as morphisms.
By  DG functor we mean a a \ka--linear functor between DG categories whose action on Hom--complexes preserves grading and commutes with differentials.
Let ${\cal A}$ be a small DG category.  We call a (left) {\cal A}--module a DG functor $M:{\cal A}\ra  \mbox{Dif}(\ka)$. We introduce a category
$\mbox{Dif}({\cal A})$ whose objects are ${\cal A}$--modules and
\[\mbox{Hom}_{\mbox{\scriptsize{Dif}}({\cal A})}(M,N):=\bigoplus_{j\in{\bf Z}} \mbox{Nat}(M,N[j])\]
where Nat stands for the set of transformations of degree $0$ between underlying graded functors (i.e. we do not assume (like in internal Hom) that transformations commute with differentials). Thus $\mbox{Dif}({\cal A})$ is a DG category [K1, Sect.~1].\newline
A natural environment for developing homological algebra in $\mbox{Dif}({\cal A})$ is its derived category ${\cal DA}$.
Let ${\cal CA}$ be the the category with the same objects as $\mbox{Dif}({\cal A})$ but
\[\mbox{Hom}_{{\cal CA}}(M,N):= Z^0(\mbox{Hom}_{\mbox{\scriptsize{Dif}}({\cal A})}(M,N))\]
i.e. we consider only morphisms of degree 0 commuting with differentials.
Then the quickest way of defining
${\cal DA}$ is just by saying that it is the localization of ${\cal CA}$ with respect to the class of quasi--isomorphisms. However, in order to get a
more concrete description of  ${\cal DA}$ it is convenient to use the formalism of Quillen model categories.
There are two Quillen model structures on  $\mbox{Dif}({\cal A})$: the projective and  injective one, in both structures the class of weak
equivalences is the class of quasi--isomorphisms. We focus here on the projective structure. In this structure every object is cofibrant while
the fibrant objects are those satisfying ``property P'' [K1, Sect.~3], [K2, Sect.~3.2]. Of course, like in any model category, every object is quasi--isomorphic (=weakly equivalent)
to a fibrant one. However, in our situation this can be made functorially. Namely, let
Let ${\cal HA}$ stands for the category with the same objects as $\mbox{Dif}({\cal A})$ but
\[\mbox{Hom}_{{\cal HA}}(M,N):= H^0(\mbox{Hom}_{\mbox{\scriptsize{Dif}}({\cal A})}(M,N))\]
i.e. this time we consider  morphisms of degree 0 commuting with differentials modulo chain homotopy. Let ${\cal HA}_p$ be the full subcategory of ${\cal HA}$
consisting of fibrant objects. Then for any $M\in\mbox{Dif}({\cal A})$ we can choose a quasi--isomorphic fibrant ${\cal A}$--module $p(M)$ in such a way
that we get a functor $p: {\cal HA}\ra {\cal HA}_p$  which is left adjoint to the forgetful functor [K2, Prop..~3.1]. Now we can describe ${\cal DA}$ more explicitly,
since the natural projection induces an equivalence of triangulated categories:
\[{\cal HA}_p\simeq  {\cal DA}.\]
In fact for practical computations in  ${\cal DA}$ the following basic properties of representable and fibrant  ${\cal A}$--modules are usually sufficient:
  \begin{fact}
Let  $h^A\in\mbox{Dif}({\cal A})$ for $A\in{\cal A}$ be the ${\cal A}$--module represented by $A$ i.e $h^A(A'):=\mbox{Hom}_{{\cal A}}(A,A')$. Then
\begin{enumerate}
\item If $M\in\mbox{Dif}({\cal A})$ is fibrant then
for any $N\in\mbox{Dif}({\cal A})$.
\[\mbox{Hom}_{{\cal DA}}(M,N)\simeq\mbox{Hom}_{{\cal HA}}(M,N)\]
\item $h^A$ is fibrant.
\item For any $N\in\mbox{Dif}({\cal A})$
\[\mbox{Hom}_{{\cal DA}}(h^A,G)\simeq H^0(G(A)).\]
\end{enumerate}
\end{fact}
The formalism of Quillen model categories also allows to construct derived functors. We shall need this construction in a very special case of ``standard functors''
produced by bimodules [K1, Sect.~6]. Let ${\cal B}$ be another small DG category.
Then an ${\cal A}$--${\cal B}$ bimodule $X$ is an
object of the category $\mbox{Dif}({\cal A}\ot{\cal B}^{op})$. With this data one can associate a pair of adjoint functors $H_X:\mbox{Dif}({\cal A})\ra
\mbox{Dif}({\cal B})$,
$T_X:\mbox{Dif}({\cal B})\ra\mbox{Dif}({\cal A})$ given by the formulae:
\[H_X(M)(B):=\mbox{{\cal A}}(X(-,B),M),\]
\[T_X(N)(A):=\mbox{coker}(\bigoplus_{B,C\in{\cal B}} N(C)\ot \mbox{Hom}_{{\cal B}}(B,C)\otimes X(A,B)\stackrel{\nu}{\ra} \bigoplus_{B\in{\cal B}}
N(B)\ot X(A,B)),\]
where $\nu(n\ot f\ot x)=N(n)(f)\ot x-n\ot X(A,f)(x)$.\newline
Now if $X(-,B)$ is fibrant object in $\mbox{Dif}({\cal A})$ for any $B\in {\cal B}$ then $T_X$ preserves fibrant objects while $H_X$ preserves quasi--isomorphisms.
Thus  the functors ${\bf R}H_X:{\cal DA}\ra{\cal DB}$,
${\bf L}T_X:{\cal DB}\ra{\cal DA}$ given by the formulae:
\[{\bf R}H_X(M):= H_X(M),\]
\[{\bf L}T_X(N):=T_X(p(N))\]
form the pair of adjoint functors between triangulated categories.\newline
We will use this formalism in Sections 4,5 in two special cases. The first is when we have a quasi--isomorphism $\phi: {\cal B}\ra{\cal A}$
(i.e. $\phi$ is a DG functor which induces an equivalence  of the cohomology categories).
Then by applying this machinery to the bimodule $X(A,B):=\mbox{Hom}_{{\cal A}}(A,\phi(B))$ we get
\begin{fact}
Let  $\phi: {\cal B}\ra {\cal A}$ be a quasi--isomorphism of small DG categories. Then the functor ${\bf R}\phi^*:{\cal DA}\ra{\cal DB}$
given by the formula
\[{\bf R}\phi^*(M)(B):=M(\phi(B))\]
is an equivalence of triangulated categories.\newline
Moreover, its inverse $({\bf R}\phi^*)^{-1}$ satisfies the property
\[({\bf R}\phi^*)^{-1}(h^B)\simeq h^{\phi(B)}\]
for any $B\in{\cal B}$.
\end{fact}
Another instance of this construction will be crucial in Section 5.
Let $X$ be an ${\cal A}$--${\cal B}$ bimodule. We introduce the DG category ${\cal B}_X$.
Its objects are those of ${\cal B}$ but
\[\mbox{Hom}_{{\cal B}_X}(B,B'):=\mbox{Hom}_{\mbox{\scriptsize{Dif}}(A)}(X(-,B'),X(-,B)).\]
This is a \ka--linear DG category and the action of $X$ on morphisms defines a functor $\alpha: {\cal B}\ra{\cal B}_X$. Now we see that
 $X$ is an ${\cal A}$--${\cal B}_X$ bimodule. Therefore if
$X(-,B)$ is fibrant for any $B\in{\cal B}$,
we have the pair of adjoint functors ${{\bf L}T_X,{\bf R}H_X}$ between the derived categories
${\cal DB}_X$ and ${\cal DA}$. The special feature of this bimodule is the following
\begin{prop}
Let $X$ be an ${\cal A}$--${\cal B}$ bimodule regarded as ${\cal A}$--${\cal B}_X$ bimodule.
Assume that for any $B\in{\cal B}$ $X(-,B)$ is fibrant and small ${\cal A}$ module.
Then the unit map $\mu: Id_{{\cal DB}_X}\ra {\bf R}H_X\circ {\bf L}T_X$ is an isomorphism.
\end{prop}
{\bf Proof: }
Observe that, since $X(-,B)$ is small  for all $B\in{\cal B}$, ${\bf R}H_X$ commutes with infinite sums, hence so does $ {\bf R}H_X\circ {\bf L}T_X$.
Thus it suffices to show that  $\mu(h^B)$ is an isomorphism. Moreover, since $h^B$ is fibrant, we have  ${\bf L}T_X(h^B)=T_X(h^B)=X(-,B)$.
Then we obtain
\[{\bf R}H_X(X(-,B))(B')=H_X(X(-,B))(B')=\mbox{Hom}_{\mbox{\scriptsize{Dif}}(A)}(X(-,B),X(-,B'))=
\]\[=\mbox{Hom}_{{\cal B}_X^{op}}(B,B')=h^B(B').\]\qed
Now we would like to apply this machinery to our graded category ${\cal P}_d^{af}$. It is natural to introduce a DG category $\mbox{Dif}^{f+}(\Gamma^d({\cal V}_A))$ which consists of DG functors from $\Gamma^d({\cal V}_A)$ (regarded as a DG category with trivial differentials) to the category $\mbox{Dif}^{f+}(\ka)$ of bounded below complexes of finite dimensional
vector spaces over \ka.  The reason for which we keep our vector spaces finite dimensional is that we want the Kuhn duality to be an equivalence. Then, as we will see, we have to restrict to bounded below complexes, since only for such complexes fibrant resolutions remain finite dimensional. Now in order to make sure that the
category  $\mbox{Dif}^{f+}(\Gamma^d({\cal V}_A))$ still has the (projective) Quillen model structure we need the fact that fibrant resolutions exist inside this
subcategory of  $\mbox{Dif}(\Gamma^d({\cal V}_A))$.
\begin{theo}
Let $F\in\mbox{Dif}^{f+}(\Gamma^d({\cal V}_A))$. Then a fibrant resolution $p(F)$ can be chosen so that  $p(F)\in\mbox{Dif}^{f+}(\Gamma^d({\cal V}_A))$.
\end{theo}
Before we start the proof let us explain the reason for which this theorem holds, since at first sight it is strange that we have bounded below projective
resolutions. To this end let us look at the simplest case of $d=1$. In this case by Prop.~2.5, ${\cal P}_1^{af}$ is equivalent to the category
of finitely generated bounded below graded $A$--modules. Now the projective periodic resolution of the trivial $A$--module \ka\ can be written as $P_{\bullet}$ for
$\bullet\leq 0$ where $P_{-2j}=A[-2pj]$, $P_{-(2j+1)}=A[-(2pj+2)]$ and $d_{-2j}=\cdot x$, $d_{-(2j+1)}=\cdot x^{p-1}$. Hence we see that, essentially because
$A$ is connected and positively graded, $P_{\bullet}$ as DG module is bounded below.
 The same phenomenon  occurs with the bar resolution which can be applied in a more general situation. Another important point which is crucial for
 extending this result to $d>1$ is that ${\cal P}_d$ has finite homological dimension [To]. These observations will guide the proof.\newline
{\bf Proof: } Since $p(F)$ has a filtration with associated object isomorphic to $H^*(F)$, it suffices to find a resolution of $H^*(F)$ inside
$\mbox{Dif}^{f+}(\Gamma^d({\cal V}_A))$. It will be more transparent to work  with a graded $S_{d,d}^{af}$--module $M:=F(A^d)$ instead of a functor $F$ which is equivalent by Proposition 2.5. Then $M$ is generated by a countable set ${c_j}_{j\geq 0}$ such that we have only a finite set $c_j$'s in each degree.
Let $M_j:=< c_k>_{k\leq j}$. Then $M=\bigcup_j M_j$ and in each degree the filtration $M_j$ stabilizes. Thus it suffices to find  a resolution
for each $M_j$ separately. Each $M_j$ is totally finite dimensional hence it has a finite filtration with subquotients concentrated in a single degree. Thus by
the horseshoe lemma it satisfies to find resolutions for all modules concentrated in a single degree.
Let $N$ be such a module (concentrated in degree 0 to fix attention). Since any positively graded element of $S_{d,d}^{af}$ acts on $N$ trivially, N is
of the form $r^*(N')$ where $N'$ is the underlying $S_{d,d}$--module and $r: S_{d,d}^{af}\ra S_{d,d}$ is the projection onto the degree 0 part.
Since $N'$ has a finite resolution by sums of tensor powers of $S_{d',d'}$--modules $\Gamma^{d'}(\ka^{d'})$, we are left with a problem of finding
a bounded below resolution of the graded $S_{d,d}^{af}$--module $\Gamma^{d}(\ka^d)$ by finite dimensional projective graded $S_{d,d}^{af}$--modules.
To this end we take the normalized bar resolution $P_{\bullet}$ of $S_{d,d}^{af}$--module $\Gamma^d({\bf k}^d)$ associated to the ring extension
$S_{d,d}\subset S_{d,d}^{af}$. Then $P_{-j}=S^{af}_{d,d}\ot(\widetilde{S^{af}_{d,d}})^{\ot j}$
where $\widetilde{S^{af}_{d,d}}:=\mbox{ker}(r)$.
Now, since $\widetilde{S^{af}_{d,d}}$ starts in degree $2$, $P_{-j}$ starts in degree $2j$. Hence $P_{\bullet}$ regarded as a DG module
is bounded below.\qed
This theorem allows us to equip $\mbox{Dif}^{f+}(\Gamma^d({\cal V}_A))$ with the projective model structure and apply to it all the constructions described earlier in this section.
In particular we will heavily use the triangulated category ${\cal D}^{f+}\Gamma^d({\cal V}_A)$
obtained from ${\cal C}^{f+}(\Gamma^d({\cal V}_A))$ by inverting quasi--isomorphisms.
We shall denote this category as ${\cal DP}_d^{af}$, for it  should be thought of as
the derived category of ${\cal P}_d^{af}$. However, in the next sections when constructing derived functors on  ${\cal DP}_d^{af}$ we will have to check carefully that they preserve our extra finiteness and boundedness below assumptions.\newline We finish this section by remarking that the fact that
$z^*$ extends to the graded functor on ${\cal P}^{f+}$ allows to further extend it to the functor (denoted by the same symbol)
\[z^*:\mbox{Dif}^{f+}(\Gamma^d{\cal V})\ra \mbox{Dif}^{f+}(\Gamma^d{\cal V}_A).\] Thus finally we get the exact functor
\[z^*:{\cal DP}_d\ra {\cal DP}_d^{af}.\]
Analogously we obtain the exact functor
\[t^*:{\cal DP}_d^{af}\ra {\cal DP}_d\]
and they are still adjoint.

\section{Formality}
In this section we introduce and study certain DG category denoted by $\Gamma^d{\cal V}_X$ whose cohomology category is $\Gamma^d{\cal V}_A$. The main result of this section, Theorem 4.2
says that this category is formal
i.e. quasi--isomorphic to its cohomology category $\Gamma^d{\cal V}_A$. This  gives rise to a derived equivalence of respective functor categories (Corollary 4.3). This derived equivalence
will be used in the next section to lifting the action of $\Gamma^d{\cal V}_A$ from cohomology to complexes which is the crucial part of the construction
of the affine derived right Kan extension.\newline
Let ${\cal P}_{dp}^d$ be the category of bifunctors of bi--degree $(dp,d)$ in the sense of [FS] (in fact ${\cal P}_{dp}^d$ is nothing but
$\mbox{Fun}(\Gamma^{pd}{\cal V}\times \Gamma^d{\cal V}^{op},{\cal V}$). We consider the bifunctor $\Gamma^d(I^{(1)}\ot I^*)\in {\cal P}_{dp}^d$  given by the formula $(W,V)\mapsto \Gamma^d(W^{(1)}\ot V^*)$.
Now we recall that ${\cal P}_{dp}^d$ like ${\cal P}_d$  has  finite homological dimension.
Thus we choose a finite projective resolution $X$ of $\Gamma^d(I^{(1)}\ot I^*))$. Since $X$ is an object of the category $\mbox{Dif}(\Gamma^{dp}{\cal V}\ot \Gamma^d{\cal V}^{op})$ we can consider the category $\Gamma^d{\cal V}_X$. We recall from Section 3 that the objects
of this category are \ka--vector spaces and
\[\mbox{Hom}_{\Gamma^d{\cal V}_X}(V,V'):=\mbox{Hom}_{{\cal P}_{dp}}(X(-,V'),X(-,V)).\]
Now it follows from the standard Ext--computations that
\begin{prop}
The assignment $V\mapsto V\ot A$ extends to an equivalence of graded categories
\[H^*(\Gamma^d{\cal V}_X)\simeq \Gamma^d{\cal V}_A.\]
\end{prop}
{\bf Proof: } Indeed, all we need is a natural in $V,V'$ calculation of Ext--groups
\[\mbox{Ext}^*_{{\cal P}_{dp}}(\Gamma^d((-)^{(1)}\ot V'^*), \Gamma^d((-)^{(1)}\ot V^*))\simeq \Gamma^d(\mbox{Hom}(V,V')\ot A).\]
This isomorphism is the ``parameterized version'' of  [FFSS,~Th.5.4] and can be easily proved by the methods used there. We leave the straightforward
details to the reader.\qed
A much deeper is the following result which is yet another incarnation of formality phenomena observed in [C4].
\begin{theo}
The assignment $V\mapsto V\ot A$ extends to an equivalence of DG categories $\phi: {\cal V}_A\simeq \Gamma^d{\cal V}_X$.
\end{theo}
{\bf Proof: } We define $\phi_{V,V'}$ as the composite of several operations.\newline
First we dualize:
\[\Gamma^d(\mbox{Hom}(V,V')\ot A))\simeq \Gamma^d(\mbox{Hom}(V'^*\ot A^*,V^*))\]
and apply the Yoneda lemma:
\[ \Gamma^d(\mbox{Hom}(V'^*\ot A^*,V^*)\simeq \mbox{Hom}_{{\cal P}_d}(\Gamma^d(-\ot V'^*\ot A^*),\Gamma^d(-\ot V^*))\]
which clearly preserves composing morphisms.\newline
Then we precompose with $I^{(1)}$
\[\mbox{Hom}_{{\cal P}_d}(\Gamma^d(-\ot V'^*\ot A^*),\Gamma^d(-\ot V^*))\ra \mbox{Hom}_{{\cal P}_{pd}}(\Gamma^d((-)^{(1)}\ot V'^*A^*),\Gamma^d((-)^{(1)}\ot V^*)) \]
which also commutes with composing.\newline
Next we lift morphisms to resolutions
\[\mbox{Hom}_{{\cal P}_{pd}}(\Gamma^d((-)^{(1)}\ot V'^*\ot A^*),\Gamma^d((-)^{(1)}\ot V)\ra
\mbox{Hom}_{{\cal P}_{pd}}(X(-,V'\ot A^*),X(-,V)).\]
The commuting of the lift with composition follows from the functoriality of $X$ (or rather its extension to the graded spaces) with respect to the
first variable. Now we recall from [C3, Prop.~3.2] that there is an element $e_d\in \mbox{Hom}_{{\cal DP}_{dp}^d}(\Gamma^d(I^{(1)}\otimes I^*),\Gamma^d(I^{(1)}\otimes I^*\ot A^*))$ with certain special properties.
Then we can realize $e_d$ as an element in $\mbox{Hom}_{{\cal P}_{dp}^d}(X(-,-),X(-,-\ot A^*))$. Hence composing with $e_d$
evaluated on $V'$ produces
for any $V,V'$ the arrow
\[\mbox{Hom}_{{\cal P}_{pd}}(X(-,V'\ot A^*),X(-,V))\ra\mbox{Hom}_{{\cal P}_{pd}}(X(-,V'),X(-,V))\]
which by the naturality of $e_d$, commutes with composition. The composite of these four arrows is a natural with respect to $V,V'$ arrow $\phi$:
\[\phi_{V,V'}:\Gamma^d(\mbox{Hom}(V,V')\ot A))\ra\mbox{Hom}_{{\cal P}_{pd}}(X(-,V'),X(-,V)).\]
The reader of [C4] will recognize  that our $\phi_{V,V'}$ is nothing but $\Phi(\Gamma^{d,V})(V')$. It was showed in the proof of [C4, Th.~3.2] that this arrow
is a quasi--isomorphism. We repeated this construction here, since for our Theorem 2.2 we need strict  commuting of our arrows with composition.\qed
Now we would like to deduce from Theorem 4.2 the derived equivalence of the respective functor categories
by using Fact 3.2. However, when we apply Fact 3.2 literally we obtain an equivalence
\[{\bf R}\phi^*:   {\cal D}\Gamma^d{\cal V}_X\simeq {\cal D}\Gamma^d{\cal V}_A.\]
Therefore we should check wether $\phi^*$ preserves f+ conditions. To this end we recall that by Fact 3.2 again, $({\bf R}\phi^*)^{-1}$ preserves representable objects. Hence, since Hom--complexes in $\Gamma^d{\cal V}^{op}_X$ are totally
finite dimensional (we use here the fact that $X$ is of finite length), $({\bf R}\phi^*)^{-1}$ preserve finite dimensionality and boundedness below
conditions. Hence we obtain
\begin{cor}
The functor
\[{\bf R}\phi^*:{\cal DP}_X:=
{\cal D}^{f+}\Gamma^d{\cal V}^{op}_X\ra {\cal D}^{f+}\Gamma^d{\cal V}_A={\cal DP}_d^{af}\]
 is an equivalence of triangulated categories.
Moreover $({\bf R}\phi^*)^{-1}(h^V)\simeq hx^V$ where $hx^V$ is the object represented by $V\in{\cal V}$ in $\mbox{Dif}(\Gamma^d{\cal V}_X)$.
\end{cor}

\section{The affine derived right Kan extension}
In this section we use the  complex $X$ regarded as $\Gamma^{pd}{\cal V}$--$\Gamma^{d}{\cal V}_X$ bimodule to construct an ``affine'' version
of the derived right Kan extension introduced in [C4].\newline
First we observe that since $X$ is totally finite dimensional, the ``standard functors'' from Section 3:
\[H_X: \mbox{Dif}(\Gamma^{pd}{\cal V})\ra \mbox{Dif}(\Gamma^d{\cal V}_X),\ \ \ \ T_X:\mbox{Dif}(\Gamma^{d}{\cal V}_X)\ra \mbox{Dif}(\Gamma^{pd}{\cal V})\]
preserve the subcategaries $\mbox{Dif}^{f+}$. Hence we pass to  the derived categories and (slightly abusing notation by using the same letters)  consider
the derived functors
\[{\bf R}H_X: {\cal DP}_{pd}\ra {\cal DP}_X\ \ \ \ {\bf L}T_X:{\cal DP}_X\ra{\cal DP}_{pd}.\]
Now we are ready for defining our ``affine'' adjunction. We define the affine precomposition with the Frobenius twist:
\[ {\bf C}^{af}:{\cal DP}_{d}^{af}\ra {\cal DP}_{d}\]
as ${\bf C}^{af}:={\bf L}T_X\circ ({\bf R}\phi^*)^{-1}$,\newline
and the affine derived right Kan extension:
\[{\bf K}^{af}:{\cal DP}_{pd}\ra {\cal DP}_d^{af}\]
as ${\bf K}^{af}:={\bf R}\phi^*\circ {\bf R}H_X$.\newline
Our terminology is justified by the following theorem (where ${\bf C}$ and ${\bf K}^r$ stand respectively for the precomposition with the Frobenius twist
and the derived right Kan extension from [C4]).
\begin{theo}
There are following isomorphisms of functors:
\begin{enumerate}
\item ${\bf C}^{af}\circ z^*\simeq {\bf C}$, $t^*\circ {\bf K}^{af}\simeq {\bf K}^r$.
\item ${\bf K}^{af}$ is right adjoint to ${\bf C}^{af}$.
\item ${\bf K}^{af}\circ{\bf C}^{af}\simeq Id_{{\cal DP}_d^{af}}$.
\end{enumerate}
\end{theo}
{\bf Proof: }
In order to get the first part we evaluate ${\bf C}^{af}\circ z^*$ on the projective generator $\Gamma^{d,V}$. We obtain
\[{\bf C}^{af}\circ z^*(\Gamma^{d,V})={\bf C}^{af}(h^V)={\bf L}T_X\circ ({\bf R}\phi^*)^{-1}(h^U)={\bf L}T_X(hx^V)=X(-,V)\simeq \Gamma^d((-)^{(1)}\ot V^*)\]
which gives the first isomorphism, since both functors commute with infinite sums. To get the second formula we evaluate $t^*\circ {\bf K}^{af}$ on any $F\in{\cal P}_{pd}$.
This time we just obtain
\[t^*\circ {\bf K}^{af}(F)=\mbox{Hom}_{{\cal P}_{pd}}(X,F)={\bf K}^r(F).\]
The second part of the theorem follows from the $\{T_X,H_X\}$ adjunction and the fact that ${\bf R}\phi^*$ is an equivalence.\newline
The third part  follows from Proposition 3.3  and again the fact that ${\bf R}\phi^*$ is an equivalence.\qed
Now the main result of [C4] follows as a formal consequence.
\begin{cor}{{\bf \ (The Collapsing Conjecture)}}\newline
For any $F,G\in{\cal P}_d$ and $i>0$, there is a natural in $F,G$ isomorphism
\[\mbox{Ext}^*_{{\cal P}_{p^id}}(F^{(i)},G^{(i)})\simeq\mbox{Ext}^*_{{\cal P}_d}(F,G_{A_i})\]
where $A_i:=A\ot A^{(1)}\ot\ldots A^{(i-1)}=\ka[x_1,\ldots,x_i]/(x_1^p,\ldots,x_i^p)$ for $|x_j|=2p^{j-1}$
and $G_{A_i}(V):=G(V\ot A_i)$.
\end{cor}
{\bf Proof: } We start with $i=1$. We get
\[\mbox{Ext}^*_{{\cal P}_{pd}}(F^{(1)},G^{(1)})=\mbox{Ext}^*_{{\cal P}_{pd}}({\bf C}(F),{\bf C}(G))
\simeq \mbox{HExt}^*_{{\cal P}_{d}}(F,{\bf K}^r{\bf C}(G))\simeq\] \[ \mbox{HExt}^*_{{\cal P}_{d}}(F,t^*{\bf K}^{af}{\bf C}^{af}z^*(G))
\simeq  \mbox{HExt}^*_{{\cal P}_{d}}(F,t^*z^*(G))= \mbox{Ext}^*_{{\cal P}_{d}}(F,G_A).\]
The case $i>1$ is obtained by iterating this computation. An important point is that $t^*z^*(G_{A_i})=G_{A_{i+1}}$ which follows
from the fact (observed already in [C1]) that the Frobenius twist
extended to the graded spaces multiplies degrees by $p$ (see also  the end of the proof of [C4, Th.~3.2]).\qed
In fact, a part of motivation for the present work was to put the Collapsing Conjecture into a more general categorical context.
The Collapsing Conjecture is essentially a statement about
the unit of the $\{{\bf C},{\bf K}^r\}$ adjunction.
Our Theorem 5.1 allows to divide its construction into two steps.
The first is  sort
of scalar extension from ${\cal P}_d$ to ${\cal P}_d^{af}$ and
we see that the unit here is the precomposition with the graded space $A$.
The second is the affine derived right Kan extension whose unit is just the identity. This point of view offers somewhat more conceptual picture of the
Collapsing Conjecture.
\section{Concluding remarks}
In this section we briefly discuss various implications and ramifications  of our work and sketch possible further developments.
As we have explained in the previous section the affine strict polynomial functors help to better understand phenomena surrounding the Collapsing Conjecture.
However, my original motivation  was more general. Namely, as I have mentioned in the Introduction, the category ${\cal P}_d^{af}$
provides conceptual explanation of various homological
phenomena in ${\cal P}_{d}$ which were observed empirically on many occasions. For example, in many calculations of the groups $\mbox{Ext}^*_{{\cal P}_{pd}}(F^{(1)},G)$ some extra structure seemed to emerge. To put it simply: all these
computations were given in terms of the graded space $A$. What is important, this phenomenon is not restricted to the case of $G=G'^{(1)}$
which is covered by the Collapsing Conjecture but  can also be observed e.g. in $\mbox{Ext}^*_{{\cal P}_{pd}}(W_{\mu}^{(1)},S_{\lambda})$
(the Ext groups between twisted Weyl and Schur functors)
  whose computing is crucial for understanding the structure of ${\cal DP}_{pd}$. This was already apparent in [C3] where the groups $\mbox{Ext}^*_{{\cal P}_{pd}}(W_{\mu}^{(1)},S_{\lambda})$ were computed in certain special case, but it will be seen much more vividly in [C5] where the case
of ``p--quotient consisting of several diagrams'' is considered. With our factorization ${\bf K}^r=t^*\circ{\bf K}^{af}$ this becomes quite natural since
we have
\[\mbox{Ext}^*_{{\cal P}_{pd}}(F^{(1)},G)\simeq\mbox{HExt}^*_{{\cal P}_{d}}(F,{\bf K}^r(G))\simeq
\mbox{HExt}^*_{{\cal P}_{d}}(F,t^*{\bf K}^{af}(G))\]
and we recall that $t^*{\bf K}^{af}(G)(V):={\bf K}^{af}(G)(V\ot A)$. In the case of Ext--groups between Weyl and Schur functors it becomes a central point, since it was shown in [C4, Prop.~4.1] that the problem of computing  $\mbox{Ext}^*_{{\cal P}_{pd}}(W_{\mu}^{(1)},S_{\lambda})$  essentially reduces to
that of finding of ${\bf K}^r(S_{\lambda})$. We will see in [C5] that indeed in practice we describe ${\bf K}^r(S_{\lambda})$ as $t^*{\bf K}^{af}(S_{\lambda})$ for certain  explicitly described affine functor ${\bf K}^{af}(S_{\lambda})$.\newline
Another mysterious fact emerging from Ext--computations for Weyl and Schur functors was that they were governed by the combinatorics of p--quotients
of Young diagrams. Here too, the category ${\cal P}_d^{af}$ provides sort of heuristic explanation. It is best seen in terms of representations
instead of functors. Namely, since $\mbox{End}_A(A^n)$ acts on $F(A^n)$ for any $F\in {\cal P}_d^{af}$, affine strict polynomial functors
produce graded (polynomial) representations of the graded algebraic group $GL_n(A)$. Since    $A$ is the group algebra for the cyclic group
${\bf Z}/p$, our group is closely related to the group $GL_n$ with coefficients in the group algebra for the infinite cyclic group.
But this group: $GL_n(\ka[x,x^{-1}])$ is nothing but the group of algebraic loops on $GL_n(\ka)$. The latter group is the affine Kac--Moody
group of type $A_n$. In particular, its simple representations are labeled by n--tuples of Young diagrams. This analogy explains why we call
our functors affine and also suggest that the combinatorics of tuples of Young diagrams should somehow organize the structure of ${\cal P}_d^{af}$. To be more precise: it seems that the relevant combinatorial structure is the set of p--tuples of Young diagrams with the total weight d.  One possible way of incorporating combinatorics into a structure of our category would be by showing that it is an ``A-highest weight category'' in the sense of [Kle] but we postpone a more systematic study of the structure of ${\cal P}_d^{af}$ to a future work.
At the time being we can only announce an explicit construction (which will be described in detail in [C5]) which explains how combinatorics
of p-tuples of Young diagrams appears in Ext--computations. Namely there exists in the category ${\cal P}_d^{af}$ a construction analogous
to that of Schur functor. It  associates to a p--tuple of Young diagrams $\{\lambda^0,\ldots,\lambda^{p-1}\}$ the affine functor
$S^{af}_{\{\lambda^0,\ldots,\lambda^{p-1}\}}\in {\cal P}_d^{af}$ by means of certain symmetrizations, antysymmetrizations etc. Then, it will
be shown in [C5] that if p--core of $\lambda$  is trivial then ${\bf K}^{af}(S_{\lambda})=S^{af}_{q(\lambda)}$ where $q(\lambda)$ is  p--quotient of $\lambda$ (for the definition and basic properties of p--cores and p--quotient consult e.g. [JK, pp. 75--76]).
At last we remark that if p--quotient of $\lambda$ consists of a single diagram $\tau$ then
\[S_{q(\lambda)}^{af}(V\ot A)\simeq  S_{\tau}(V)\]
up to shift. This will allow to see results of [C3] as special case of those of [C5].\newline

Now I would like to make some comments on choices I have made and certain limitations of the present work.
In fact there are several different reasonable candidates for ${\cal V}_A$. In a sense a more natural choice would be the category with objects of the form
$V^{\bullet}\ot A$ where $V^{\bullet}$ is a graded \ka--space or even the whole category of graded free finitely generated $A$--modules. In each case it is possible to
develop a parallel  theory of affine strict polynomial functors enjoying some advantages over the variant we have chosen. For example admitting
graded vector spaces would make discussion of extending functors to the graded spaces more transparent: we would develop the whole theory
including affine derived Kan extension for functors on graded spaces and only at the very end we would compare the graded and ungraded versions
of functor categories end deduce the ungraded theory from the graded one. Also, allowing all graded free finitely generated $A$--modules would
make the definition of Kuhn duality more straightforward. Despite all  advantages of these alternative approaches my feeling is that it is best
to keep the category ${\cal V}_A$ as small as possible (as long as it is able to capture the affine derived right Kan extension). An additional reason for not allowing all graded free finitely generated $A$--modules is that $H^*(X)$ is not a functor on that category. Of course this is not a serious problem since it is defined on the equivalent full subcategory consisting of those of the form $V^{\bullet}\ot A$ but this suggests that
this larger category is not a natural environment for our constructions. Therefore I have decided to keep our category  smallest but it is possible that in some future applications the alternative approaches will turn out to better suited.\newline
Also, in contrast to [C4] I decided not to consider multiple Frobenius twists. The reason was mainly just not to overload notation by adding another index $i$ while usually in applications like Cor.~5.2 the case of multiple twists can be deduced from the case of a single twist just by iteration. \newline Yet another difference with [C4] is of much more fundamental nature. Namely we do not consider the affine derived left Kan extension. Of course we can define ${\bf C}^{af}_l:=(-)^{\#}\circ{\bf C}^{af}\circ (-)^{\#}$, ${\bf K}^{af}_l:=(-)^{\#}\circ{\bf K}^{af}\circ (-)^{\#}$ and ${\bf K}^{af}_l$ is left adjoint to ${\bf C}^{af}_l$.
But the problem is that ${\bf C}^{af}_l$ does not preserve $\mbox{Dif}^{f+}$ (its natural domain is $\mbox{Dif}^{f-}$). On the other hand if we drop finitness/boundedness assumptions the Kuhn
duality is not an equivalence anymore. Moreover, which is even more discouraging, even for bounded complexes ${\bf C}^{af}$ does not
commute with the Kuhn duality. To see this we recall from Prop.~2.4.3 that $(\chi_0)^{\#}=\chi_0[-2(p-1)]$. But it is easy to see that ${\bf C}^{af}(\chi_0)$ is not bounded above, hence its dual is not bounded below. For these reasons one can doubt whether ${\bf C}^{af}$ has a left adjoint in any reasonable context.


\end{document}